\documentclass[english,russian]{amsart}
\usepackage{amssymb,amsmath}
\newtheorem{thm}{Theorem}
\newtheorem{crl}[thm]{Corollary}
\newtheorem{prp}[thm]{Proposition}
\newtheorem{lm}[thm]{Lemma}

\newcommand{\der }{\partial}

\newcommand{\sudda}[1]{}

\begin{document}

\title{$2p$-commutator on differential operators of order $p$}

\author{Askar Dzhumadil'daev}
\address
{Kazakh-British University, Tole bi 59, Almaty,050000,Kazakhstan}
\email{dzhuma@homail.com} 

\subjclass{16S32}

\keywords{Weyl algebra, differential operator, N-commutator, polynomial identity, Amitzur-Levitzki identity}

\maketitle

\begin{abstract} We show that a space of one variable differential operators of order $p$ admits non-trivial $2p$-commutator and the number $2p$ here can not be improved. 
\end{abstract}

Let $A$ be an associative algebra over a field $K$ of characteristic $0.$ Let $f=f(t_1,\ldots,t_n)$ be some non-commutative associative polynomial. Say that $f=0$ is {\it an identity} of $A$ if $f(a_1,\ldots,a_n)=0$ for any substitutions $t_i:=a_i\in A.$ Let $s_n$ be a skew-symmetric associative non-commutative polynomial 
$$s_n(t_1,\ldots,t_n)=\sum_{\sigma\in Sym_n} sign\,\sigma\,t_{\sigma(1)}\cdots t_{\sigma(n)}.$$
For example,
$$s_2(t_1,t_2)=t_1t_2-t_2t_1=[t_1,t_2]$$
is a Lie commutator. 

Suppose that an associative commutative algebra  $U$ has  $n$ commuting derivations $\der_1,\ldots\der_k.$ A linear span of linear operators of a form 
$u\der_{i_1}\ldots\der_{i_p},$ where $1\le i_1,\ldots,i_p\le k,$ is denoted $D_k^{(p)}(U).$ Let 
$D_k(U)=\cup_{p\ge 0} D_k^{(p)}(U)$ be space of differential operators on $U$ generated by derivations $\der_1,\ldots,\der_k.$  In case of $k=1$ we reduce notation $\der_1$ to $\der.$ 

It is known that  $D_k(U)$ can be endowed by a structure of associative algebra. A multiplication of the algebra $D(U)$ is given as a composition of differential operators. For example, if $k=1,$ then 
$$u\der^p \cdot  v\der^l=\sum_{s=0}^p{p\choose s} u\der^s(v)\der^{p+l-s}.$$
Certainly this construction can be easily generalized for algebras with several derivations. 

We can consider   $D_k^{(p)}(U)$ as a space of differential operators of order $p.$
Well known, that any differential operator of first order is a derivation and a space of derivations $Der(U)=D_k^{(1)}(U)$ 
forms Lie algebra under commutator,
$$u\der_i, v\der_j\in D_k^{(1)}(U)\Rightarrow s_2(u\der_i,v\der_j)=u\der_i \cdot v\der_j -v\der_j\cdot u\der_i\Rightarrow$$
$$ s_2(u\der_i,v\der_j)
=u\der_i(v)\der_j-v\der_j(u)\der_i\in D_k^{(1)}(U).$$

Main example of the algebra of differential operators appears in the case $U=K[x_1,\ldots,x_k]$ 
and $\der_i={\der}/{\der x_i},$ $i=1,\ldots,k,$  are partial differential operators. Recall that 
action of $\der_i$ on a monom $x^{\alpha}=x_1^{\alpha_1}\cdots x_k^{\alpha_k},$ where $\alpha=(\alpha_1,\ldots,\alpha_k)\in {\bf Z}_0^k,$ is defined by 
$$\der_i x^{\alpha}=\alpha_i x^{\alpha-\epsilon_i}.$$
Here ${\bf Z}_0$ is a set of non-negative integers and $\epsilon_i=(0,\ldots,0,1,0,\ldots,0)\in {\bf Z}_0^k$ (all components  of  $\epsilon$ except $i$-th are $0$). 

Denote by $A_k$ an algebra of differential operators on polynomials algebra $K[x_1,\ldots,x_k]$ 
generated by $k$ commuting derivations $\der_1,\ldots,\der_k.$ 
The algebra $A_k$ is called {\it Weyl} algebra.  Let $A_k^{(p)}=\langle u\der^{\alpha} | |\alpha|=p\rangle$ be subspace of $A_k$ consisting differential operators of $p$-th order. 

Let us consider $A_k^{(p)}$ as $N$-ary algebra  under $N$-ary multiplication $s_N,$ 
$$s_N(X_1,\ldots,X_N)=\sum_{\sigma\in Sym_N} sign\,\sigma\, X_{\sigma(1)}\cdots  X_{\sigma(N)}.$$
In general this notion is not correct. Might happen that $s_N$ is not well-defined on $A_k^{(p)},$
$$s_N(X_1,\ldots,X_n)\notin A_k^{(p)}$$
for some $X_1,\ldots,X_N\in A_k^{(p)}.$  We say that $A_k^{(p)}$ {\it admits $N$-commutator} $s_N,$ if 
$$s_N(X_1,\ldots,X_N)\in A_k^{(p)}$$
for any $X_1,\ldots,X_N\in A_k^{(p)}.$ 

In \cite{Dzh} it was proved that the space of differential operators of first order  $A_n^{(1)}$ in addition to   Lie commutator $s_2$ 
admits {$(n^2+2n-2)$-commutator} and that $s_N=0$ is identity if $N\ge n^2+2n.$
Let $Mat_n$ be an algebra of $n\times n$ matrices. Amitzur-Levitzky theorem states that $Mat_n$ satisfies the identity $s_{2n}=0$ and it is a minimal identity \cite{Amitzur-Levitzky}.  Note that Weyl algebra has no polynomial identity except associativity. So, to construct non-trivial identities we have to consider smaller subspaces of Weyl algebra. 

The aim of our paper is to establish that the space of one variable differential operators of order $p$ admits $2p$-commutator. 
The number $2p$ here can not be improved: if $N>2p,$ then $s_{N}=0$ is identity on $A_1^{(p)}$;
if $N<2p,$ then $s_{N}$ is not well-defined on $A_1^{(p)}$; if $N=2p,$ then $s_N$ is well-defined on $A_1^{(p)}$ and non-trivial. Obtained $2p$-ary algebra $A_1^{(p)}$ under multiplication $s_{2p}$ is simple and left-commutative. In particular, the $2p$-algebra $(A_1^{(p)},s_{2p})$ is  homotopical $2p$-Lie. 
To formulate exact result we have to introduce some definitions. 

Let us given an $n$-ary algebra $(A,\psi)$  with $n$-ary skew-symmetric multiplication $\psi:\wedge^n A\rightarrow A.$ Say that $A$ has  {\it $(2n-2,1)$-type identity}  (in 
\cite{DzhVronskian} it is called $(n-1)$-left commutative) if it satisfies the identity
$$\sum_{\sigma\in S^{(2n-2,1)}}sign\,\sigma\, \psi(a_{\sigma(1)},\ldots,a_{\sigma(n-1)},\psi(a_{\sigma(n)},\ldots,a_{\sigma(2n-2)},a_{2n-1}))=0$$
Say that $(A,\omega)$ satisfies  {\it $(1,2n-2)$-type identity}, if 
$$\sum_{\sigma\in S^{(1,2n-2)}}sign\,\sigma\, \psi(a_{1}, a_{\sigma(2)},\ldots,a_{\sigma(n-1)},\psi(a_{\sigma(n)},\ldots,a_{\sigma(2n-1)}))=0,$$
for any $a_1,\ldots,a_{2n-1}\in A.$ Here 
$$S^{(2n-1,1)}=\{\sigma\in S_{n-1,n} | \sigma(2n-1)=2n-1\},$$
$$S^{(1,2n-1)}=\sigma\in S_{n-1,n} | \sigma(1)=1\},$$
where 
$$S_{n-1,n}=\{\sigma\in S_{2n-1} | \sigma(1)<\cdots \sigma(n-1), \; \sigma(n)<\cdots<\sigma(2n-1)\}$$
is a set of shuffle $(n-1,n)$-permutations on a set $\{1,2,\ldots,2n-1\}.$ Call $n$-algebra $(A,\psi)$ 
{\it left-commutative} if it satisfies the $(2n-2,1)$-type identity. Similarly, it is called {\it right-commutative} if it has the $(1,2n-2)$-type identity. In fact, these two notions are equivalent (Lemma \ref{lcom and rcom}). 

Say that $(A,\psi)$ is {\it homotopical $n$-Lie} \cite{hanlon} if it satisfies the following identity
$$\sum_{\sigma\in S_{n-1,n}}sign\,\sigma\, \psi(a_{\sigma(1)},\ldots,a_{\sigma(n-1)},\psi(a_{\sigma(n)},\ldots,a_{\sigma(2n-1)}))=0.$$

For  $k$-ary algebra $(A,\psi)$ with $k$-multiplication $\psi:\wedge^k A\rightarrow A$ and for a subspace $I\subseteq A$ say that $I$ is {\it ideal} of $A,$ if $\psi(a_1,\ldots,a_{k-1},b)\in I,$ for any $a_1,\ldots,a_{k-1}\in A, b\in I.$ Say that $A$ is {\it simple,} if it has no ideal except $0$ and $A.$ 

In our paper we prove the following result.
  
\begin{thm} \label{main} Let $A_1=D(K[x])$ be one variable  Weyl algebra  over a field $K$ of characteristic $0.$  Then
\begin{itemize}
\item $s_{2p+1}=0$ is a polynomial identity on $A_1^{(p)}.$  
\item any polynomial identity of degree no more than $2p$ follows from the associativity one
\item $s_N$ is not well-defined on $A_1^{(p)}$ if $N<2p$
\item $s_{2p}$ is well-defined  and non-trivial  operation on  $A_1^{(p)}$
\item for any $u_1,\ldots,u_{2p}\in K[x],$ the following formula holds 

$$s_{2p}(u_1\der^p,\cdots,u_{2p}\der^p)=\lambda_p
\left|\begin{array}{cccc}u_1&u_2&\cdots&u_{2p}\\
\der(u_1)&\der(u_2)&\cdots&\der(u_{2p})\\
\vdots&\vdots&\cdots&\vdots\\
\der^{2p-1}(u_1)&\der^{2p-1}(u_2)&\cdots&\der^{2p-1}(u_{2p})\\
\end{array}\right|\der^p,$$

\bigskip

\noindent where $\lambda_p$ is a positive integer
\item the $2p$-algebra $(A_1^{(p)}, s_{2p})$ is simple and left-commutative.
\end{itemize}
\end{thm}

\begin{crl} If $k>2p,$ then  $s_k=0$ is a polynomial identity on $A_1^{(p)}.$ 
\end{crl}

\begin{crl} The $2p$-algebra $(A_1^{(p)}, s_{2p})$ is right-commutative 
\end{crl}

{\bf Proof.} It follows from Lemma \ref{lcom and rcom}. 

\begin{crl} The $2p$-algebra $(A_1^{(p)}, s_{2p})$ is homotopical $2p$-Lie. 
\end{crl}

{\bf Proof.} By Corollary 2.2 of \cite{DzhVronskian} the algebra $(A_1^{(p)},s_{2p})$ is homotopical $n$-Lie. 

\begin{crl} Any polynomial identity of  Weyl algebra $A_n$ follows from the associativity identity. 
\end{crl}

This result is known. For example it  follows from \cite{Martindale}. 

{\bf Proof.} Suppose that $A_n$ has some polynomial identity $g=0$  that does not follow from associativity identity. We can assume that $g$ is multilinear. Suppose that it has 
degree $deg\,g=d.$ Then $g=0$ induces a polynomial identity for any subspace of $A_n$. 
For example $g=0$ should be identity for $A_1^{(p)}.$ Take $p$ such that $2p>d.$ 
We obtain contradiction with the minimality of identity $s_{2p}=0$  for $A_1^{(p)}.$ 

\begin{crl} Let $U$ be an associative commutative algebra with a derivation $\der.$ Then 
 $s_{2p}$ is a  $2p$-commutator of $D^{(p)}(U)$ and  $s_N=0$ is identity on $D^{(p)}(U)$ for any $N>2p.$
\end{crl}

Proof of theorem \ref{main}  is based on super-Lagrangians calculus. We  do in next section.

\section{Super-Lagrangians algebra} 

Let  ${\bf Z}_0$  be set of non-negative integers, 
 $E$  set of sequences with non-negative integer components, and
$$E_{k}=\{\alpha=(\alpha_1,\alpha_2,\ldots,\alpha_k) | 0\le \alpha_1<\alpha_2<\cdots<\alpha_k, \; 
\alpha_i\in {\bf Z}_0\},$$
$$E_{k,0}=\{\alpha\in E_k |  \alpha_1=0\},$$
$$E_{k}(l)=\{\alpha\in E_k | |\alpha|=\sum_{i=1}^k \alpha_k=l\},$$
$$E_{k,0}(l)=\{\alpha\in E_{k,0} | |\alpha|=\sum_{i=1}^k \alpha_k=l\}.$$
We endow $E_k$ by lexicographic order, $\alpha\le \beta$ if $\alpha_1=\beta_1,\ldots,\alpha_{i-1}=\beta_{i-1},$ but $\alpha_i<\beta_i.$ This order is prolonged to order on $E$ by $\alpha<\beta$ if $\alpha\in E_k, \beta\in E_l, k<l.$

Let us consider Grassman algebra $\mathcal U$  generated by 
formal symbols $\der^{i}(a),$ where $i\in {\bf Z}_0.$ We suppose that the  generator 
$a$ is odd  and the derivation $\der$ is even. So, elements $\der^i(a)$ are odd for any $i\in{\bf Z}_0.$

For  $\alpha=(\alpha_1,\alpha_2,\ldots, \alpha_k)\in E_{k}$ set 
$$a^{\alpha}=\der^{\alpha_1}(a_1)\cdots \der^{\alpha_k}(a_k).$$
The algebra $\mathcal U$  is super-commutative and associative,
$$a^{\alpha} a^{\beta}=(-1)^{k l} a^{\beta} a^{\alpha}.$$
$$
a^{\alpha} (a^{\beta}a^{\gamma})=(a^{\alpha} a^{\beta})a^{\gamma},$$
for any $\alpha\in E_k, \beta\in E_l, \gamma\in E_s.$ In particular, $a^\alpha a^\beta=0$,
if $\alpha$ and $\beta$ have common components. 
For example, 
$$a^{(2,3,5)}a^{(1,3)}=0, \quad a^{(1,2,3,5)} a^{(0,4)}=-a^{(0,1,2,3,4,5)}.$$

Let $\mathcal L$ be an algebra of super-differential operators on $\mathcal U$ under composition. Then operators of a form $a^{\alpha} \der^{i},$ where $\alpha\in E, i\in{\bf Z}_0,$  collect a base of  $\mathcal L.$ 
Composition of operators is defined as usual 
$$u\der^k\cdot v\der^l=\sum_{i=0}^k{k\choose i} u\der^{i}(v)\der^{k+l-i},$$
where elements $u\der^i(v)\in {\mathcal U}$ are calculated in terms of super-multiplication in super-algebra $\mathcal U.$ 
 For example, if $X=a^{(2,4,5)}\der^2$ and $Y=a^{(0,1,3)}\der^3,$ then 
$$\der(a^{(0,1,3)})=a^{(0,2,3)}+a^{(0,1,4)},$$
$$\der^2(a^{(0,1,3)})=\der(\der(a^{(0,1,3)}))=\der
(a^{(0,2,3)}+a^{(0,1,4)})=a^{(1,2,3)}+a^{(0,2,4)}+a^{(0,2,4)}+a^{(0,1,5)}=$$ $$
a^{(1,2,3)}+2a^{(0,2,4)}+a^{(0,1,5)},
$$
and
$$X\cdot Y=a^{(2,4,5)}a^{(0,1,3)}\der^5+2a^{(2,4,5)}\der(a^{(0,1,3)})\der^4+a^{(2,4,5)}
\der^2(a^{(0,1,3)})\der^3=$$

$$a^{(0,1,2,3,4,5)}\der^5+2a^{(2,4,5)}(a^{(0,2,3} +a^{(0,1,4)})\der^4+a^{(2,4,5)}
(a^{(1,2,3)}+2a^{(0,2,4)}+a^{(0,1,5)})\der^3=$$

$$a^{(0,1,2,3,4,5)}\der^5,$$
since 
$$
a^{(2,4,5)}a^{(0,2,3)}= a^{(2,4,5)}a^{(0,1,4)}= a^{(2,4,5)}a^{(1,2,3)}=  a^{(2,4,5)}a^{(0,2,4)}=  a^{(2,4,5)}a^{(0,1,5)}=0.$$

Let $X=\sum_{i=k}^l X_i\in {\mathcal L},$ where 
$X_i=(\sum_{\alpha\in E} \lambda_{\alpha,i}a^{\alpha})\der^i, k\le i\le l$ 
and $X_k\ne 0.$ Take $\beta\in E$ such that $\lambda_{\beta,k}\ne 0$ and 
$\lambda_{\alpha,k}=0$ if $\alpha>\beta.$ So, $X$ has highest term 
$\lambda_{\beta,k}x^{\beta}\der^k.$ Call it {\it leader} of $X$ and denote $leader(X).$  
For example, 
$$X=2a^{(0,1,5)}\der^2+5a^{(1,2,3)}\der^3-3a^{(0,2,4)}\der^2\Rightarrow 
leader(X)=-3a^{(0,2,4)}\der^2.$$

Denote by $ {U}_k$ a linear span of base elements $a^{\alpha},$ where $\alpha\in E_k.$ 
Similarly  define linear spaces $U_{k,0}$ $U_{k}(n)$ and $U_{k,0}(n)$ as linear span 
of base elements $a^{\alpha},$ where correspondingly $\alpha\in E_{k,0},$  $\alpha\in E_{k}(n),$  and $\alpha\in E_{k,0}(n)$  

Let 
$U_k^+\subset U_k$ and  $U_k^+(n) \subset U_k(n)$ are subsets generated by linear combinations of $e^{\alpha}$ with non-negative integer coefficients,
$$U_k^+=\{\sum_{\alpha\in E_k}\lambda_{\alpha} a^{\alpha} | \lambda_\alpha\in {\bf Z}_0\},$$
$$U_k^+(n)=\{\sum_{\alpha\in E_k(n)}\lambda_{\alpha} a^{\alpha} | \lambda_\alpha\in {\bf Z}_0\}.$$
Note that $U_k^+, U_k^+(n)$ are semigroups under addition, 
$$0\in U_k^+, 0\in U_k^+(n),$$
and
$$u,v\in U_k^+\Rightarrow u+v\in U_k^+,$$
$$u,v\in U_k^+(n)\Rightarrow u+v\in U_k^+(n).$$

Let 
$$L_k=\langle a^\alpha\der^i | \alpha\in E_k, i\in {\bf Z}_0 \rangle,$$
$$L_k(n)=\langle a^{\alpha}\der^i |\; i+ |\alpha|=n, \alpha\in E_k, i\in {\bf Z}_0 \rangle.$$
Denote by ${\mathcal L}^{(\ge p)}$ a space of differential operators of order no less than $p.$

\begin{prp} \label{grading}  For any $p\ge 0$ the subspace ${\mathcal L}^{(\ge p)}$ generates left-ideal 
on the  algebra  ${\mathcal L},$ 
$${\mathcal L} {\mathcal L}^{(p)}\subseteq {\mathcal L}^{(p)}.$$
Algebras $\mathcal U$ and $\mathcal L$ are graded,
$$U_k(n) U_l(m)\subseteq U_{k+l}(n+m),$$ 
$$L_k(n) L_l(m)\subseteq L_{k+l}(n+m),$$ 
$$U_k(n) L_l(m)\subseteq L_{k+l}(n+m),$$ 
$$L_k(n) U_l(m)\subseteq L_{k+l}(n+m),$$ 
for any $k,l,n,m\in {\bf Z}_0.$ 
\end{prp}

{\bf Proof.} Evident. 

\begin{lm} \label{bir1} Let $p \ge 0.$  If $u \in { U}_k(n),$ then $a\der^p (u)\in U_{k+1,0}(n+p).$ Moreover, if $u\in U^+_k(n),$ then $a\der^p(u)\in U^+_{k+1,0}(n+p).$
\end{lm}

{\bf Proof.} Our Lemma is an easy consequence of the following statements:
$$u\in U_k(n)\Rightarrow \der(u)\in U_k(n+1),$$
$$u\in U_k^+(n)\Rightarrow \der(u)\in U_k^+(n+1).$$
To prove these statements we use induction on $p.$

For $p=0$ our statement is trivial.
Let $p=1.$
If $u=a^{\alpha}=\der^{\alpha_1}(a)\cdots \der^{\alpha_k}(a),$ then by Leibniz rule 
$\der(u)$ is a sum of monoms of a form 
$$u_i=\der^{\alpha_1}(a)\cdots \der^{\alpha_{i-1}}(a) \der^{\alpha_i+1}(a) \der^{\alpha_{i+1}}(a)\cdots \der^{\alpha_k}(a), \qquad 1\le i\le k.$$
If $\alpha_{i+1}=\alpha_i+1,$ then by super-commutativity condition  $u_i=0.$ If $\alpha_{i+1}>\alpha_i+1,$ then $u_i$ is a  base monom.    Therefore, if $\alpha\in E_k(n),$ then 
$\der(a^{\alpha})$ is a linear combination of base monoms $a^{\beta},$ where $\beta\in E_{k}(n+1)$ with coefficients that are equal to $0$ or $1.$ Hence
$$u\in U_k(n)\Rightarrow \der(u)\in U_{k}(n+1),$$
$$u\in U^+_k(n)\Rightarrow \der(u)\in U_{k}^+(n+1).$$
So, base of induction is valid. 

Suppose that 
$$u\in U_k(n)\Rightarrow \der^{p-1}(u)\in U_{k}(n+p-1).$$
Then as we established above 
$$\der^p(u)=\der( \der^{p-1}(u))\in U_{k}(n+p)$$
By similar reasons 
$$u\in U^+_k(n)\Rightarrow \der^{p-1}(u)\in U_{k}^+(n+p-1)\Rightarrow \der^p(u)=\der(\der^{p-1}(u))
\in U_{k}^+(n+p).$$

\begin{lm} \label{power} For any $k\in {\bf Z}_0$ the $k$-th power 
 $(a\der^p)^k \in\mathcal L$ is a linear combination with non-negative integer coefficients  of operators of a form $a^\alpha \der^i,$ where $\alpha\in E_k,$  $|\alpha|+i=pk$ and $i\ge p.$ 
\end{lm}

{\bf Proof.} By grading property of $\mathcal U$ and $\mathcal L$ (Proposition \ref{grading}) 
it is clear that 
$(a\der^p)^k$ is a linear combination of super-differential operators of a form $a^\alpha\der^i,$ 
where $\alpha\in E_k(pk-i)$ and $i\ge p.$  By Lemma \ref{bir1}  coefficients are non-negative integers. 

\begin{lm} \label{N>2p} If $N>2p,$ then $(a\der^p)^N=0$ and 
$$(a\der^p)^{2p}=\lambda_pa^{(0,1,2\ldots,2p-1)}\der^p,$$
for some non-negative integer $\lambda_p.$ 
\end{lm}

{\bf Proof.} If $\alpha\in E_N,$ and $N= 2p+1,$ then 
$$|\alpha|\ge \sum_{i=0}^{N-1}i=N(N-1)/2= (2p+1)p=pN.$$ 
 Therefore, by Lemma \ref{power} $(a\der)^{2p+1}=0.$ 
 So, $(a\der^p)^N=0,$  if $N>2p.$ 
 
 If $N=2p$ and $\alpha\in E_N$ then by the same reasons,
 $$|\alpha|\ge p(2p-1),$$
 and
 $$(a\der^p)^N=leader((a\der^p)^N)=\lambda_p a^{(0,1,\ldots,2p-1)}\der^p,$$
 for some $\lambda_p\in {\bf Z}_0.$ 
 
 \bigskip
 
 To prove Theorem \ref{main} we have to establish that $\lambda_p>0.$ It will be done in next section.

 \section{Positivity of $\lambda_p$}

\begin{lm}\label{maximal}
Let $\delta(k)$ be maximal element in $E_{k+1,0}(pk).$ Then 
$$\delta(k)=\left\{\begin{array}{ll}
(0,p-l,p-l+1,\ldots,p-1,\; p+1,\ldots,p+l-1, p+l), & \mbox{ if $k=2l,$}\\
(0,p-l,p-l+1,\ldots,p-1, p, p+1,\ldots,p+l-1, p+l), & \mbox{ if $k=2l+1.$}
\end{array}\right.$$
\end{lm}

{\bf Proof.} Note first of all that $\delta(k)\in E_{k+1,0}(pk).$ Indeed, 
$$|\delta(k)|=\left\{
\begin{array}{cl}
2pl, & \mbox{ if $k=2l$ is even,}\\
p(2l+1),& \mbox{ if $k=2l+1$ is odd} 
\end{array}\right. \Rightarrow |\delta(k)|=pk.
$$

Suppose that $\beta\ge \delta(k)$ for some $\beta=(\beta_1,\ldots,\beta_{k+1})\in E_{k+1,0}(pk).$ Then
$$
\beta_2\ge p-l,
$$
where $l=\lfloor n/2\rfloor.$

For $1<i\le k$ let us call $\beta_{i+1}-\beta_i$ as $i$-rise of $\beta$ and denote $r_i(\beta).$ 
If $r_i(\beta)\ge 3$ for some $1<i\le k,$ then  we can find $\gamma\in E_{k+1,0}(pk)$ such that $\beta<\gamma.$ Take for example, 
 $\gamma_j=\beta_j,$ if $j\ne i,i+1$ and $\gamma_i=\beta_i+1,\gamma_{i+1}=\beta_{i+1}-1.$ Therefore,
$$r_i(\beta)\le 2, \qquad 1<i\le k.$$

If $r_i(\beta)=2$ for some $i$ then $r_j(\beta)=1$ for any $j\ne i,$ $1<j\le k.$ Let us prove it by contradiction. Suppose that 
$r_i(\beta)=2$ and $r_j(\beta)=2$ for $i\ne j, 1<i,j\le k.$ Then there exists $\mu\in E_{k+1,0}(pk),$  such that $\beta<\mu.$ Take for example,
$\mu_s=\beta_s$, if $s\ne i,j,$ and $\mu_i=\beta_i+1,$ $\mu_{j+1}=\beta_{j+1}-1.$

Let $k=2l+1.$ If $r_{s+1}(\beta)>1,$ for some $0\le s\le k-1.$  then $\beta_{2+s}>p-l+s.$ Therefore, 
$|\beta|>\sum_{i=p-l}^{p+l}i >pk.$
Hence, $r_i(\beta)=1$ for any $1<i\le ,$ and, $\beta=\delta(k).$ 

Let $k=2l.$ If $r_{s+1}(\beta)>1,$ for some $0\le s\le l-1,$ then  
$\beta_{2+s}>p-l+s.$  Therefore,
$$\sum_{i=1}^{s+1} \beta_i\ge \sum_{i=1}^{s+1} \delta(k)_i,$$
$$\sum_{j=s+2}^{l+1}\beta_s>\sum_{j=s+2}^{l+1}\delta(k)_j,$$
$$\sum_{t=l+2}^{2l+1}\beta_t \ge \sum_{t=l+2}^{2l+1}\delta(k)_t.$$ 
Hence,  
$$|\beta|=\sum_{i=1}^{s+1} \beta_i+\sum_{j=s+2}^{l+1}\beta_j+\sum_{t=l+2}^{2l+1}\beta_t 
>\sum_{i=1}^{s+1} \delta(k)_i+\sum_{j=s+2}^{l+1}\delta(k)_j+\sum_{t=l+2}^{2l+1}\delta(k)_t =
|\delta(k)|=pk.$$
If $r_{s+1}(\beta)>1,$ for some $l< s\le k+1,$ then  
$\beta_{2+s}>p-l+s,$  and, 
 $$\sum_{i=1}^{s+1} \beta_i\ge \sum_{i=1}^{s+1} \delta(k)_i,$$
 $$\sum_{j=s+2}^{2l+1} \beta_j>\sum_{j=s+2}^{2l+1} \delta(k)_i.$$
Therefore,
$$|\beta|=\sum_{i=1}^{s+1} \beta_i+\sum_{j=s+2}^{2l+1}\beta_j>
\sum_{i=1}^{s+1} \delta(k)_i+\sum_{j=s+2}^{2l+1}\delta(k)_j=|\delta(k)|=pk.$$
Hence 
$$r_{s+1}(\beta)>1\Rightarrow s=l,$$
and $\beta=\delta(k).$
$\square$

Recall that  $\alpha=(\alpha_1,\ldots,\alpha_k)\in {\bf Z}_0^k$ is called {\it composition} of $n$ 
with length $k$ if $\sum_{i=1}^k \alpha_i=n.$ Denote by $C_k(n)$ set of compositions of $n$ of length $k.$ For $\alpha\in C_k(n)$ denote by $sort(\alpha)$ the composition $\alpha$ written  
in non-decreasing order. Note that $sort(\alpha)$ gives us a partition of $n.$ 
 For example, $sort((2,0,2,3,1))=(0,1,2,2,3).$ 
For $\sigma=(0,\sigma_2,\ldots,\sigma_{k+1})\in E_{k+1,0}(n)$ set  $\bar\sigma=(\sigma_2,\ldots,\sigma_k)\in E_{k}(n).$

For $\alpha\in E_k,\beta\in E_{k+1}$ set 
$$M(\alpha,\; \beta)=\{\gamma\in E_k | sort(\alpha+\gamma)=\bar\beta\}.$$ 

For $\alpha\in {\bf Z}_0^k, \beta\in {\bf Z}_0^l$ define $\alpha\smallsmile \beta\in {\bf Z}_0^{k+l}$ as a prepend $\alpha$ to $\beta$
$$\alpha\smallsmile \beta=(\alpha_1,\ldots,\alpha_k,\beta_1,\ldots,\beta_l).$$
Let  
$${\bf 0}_0=(\;),$$
$${\bf 0}_i=\underbrace{(0,0,\ldots,0)}_{\mbox{$i$  times}},\qquad i>0.$$
For $\alpha\in {\bf Z}_0^k$ set 
$${|\alpha|\choose \alpha}=\prod_{i=1}^k {\alpha_1+\cdots+\alpha_k\choose \alpha_1,\ldots,\alpha_k}=\frac{(\alpha_1+\cdots+\alpha_k)!}{\alpha_1!\cdots\alpha_k!}.$$

Let 
$$G_0=\{(\;)\},$$
$$G_k=\{(i)\smallsmile {\bf 0}_{i-1}\smallsmile \alpha | \alpha\in G_{k-i}, \quad i=1,2,\ldots,k\}, \qquad k>0.$$
{\bf  Example.} 
$$G_1=\{(1)\},\quad  G_2=\{(2,0),(1,1)\}, \quad G_3=\{(3,0,0),(2,0,1), (1,2,0),(1,1,1)\}.$$

\begin{lm}\label{a1}
If $k=2l-1$ is odd, 
$$M(\delta(k-1), \delta(k))=\{(p-l+i)\smallsmile {\bf 0}_{i-1}\smallsmile \alpha\smallsmile {\bf 0}_{l-1} |\alpha\in G_{l-i},
i=1,2,\ldots, l \}.$$
If $k=2l$ is even,
$$M(\delta(k-1),\delta(k))=\{(p-l)\smallsmile {\bf 0}_{l-1}\smallsmile \alpha |\alpha\in G_{l}\}.$$
\end{lm}

{\bf Proof.} Evident. 

{\bf Example.}  If $p=5,$ then 
$$M(\delta(2),\delta(3))=M((0,4,6),(0,4,5,6))=\{(4, 1,0), (5,0,0 )\},$$
$$M(\delta(3),\delta(4))=M((0,4,5,6),(0,3,4,6,7))=\{(3,0, 1,1), (3,0,2,0 )\}.$$

\begin{lm}\label{a2}
$$\sum_{\alpha\in G_k}sign\,(\alpha+(0,1,\ldots,k-1)){k\choose \alpha}=1$$
\end{lm}

{\bf Proof.} Induction on $k.$ For $k=1$ our statement is evident. Suppose that it is true for $k-1.$ 
Note that
$$G_k=\cup_{i=1}^k  \{(i)\smallsmile {\bf 0}_{i-1}\smallsmile G_{k-i}\}.$$
For $\alpha\in G_{k-i},$ 
$$(i)\smallsmile {0}_{i-1}\smallsmile \alpha+(0,1,\ldots,k-1)=(i,1,2,\ldots, i-1,\alpha_1+i,\ldots,\alpha_{k-i}+k-1),$$
and,
$$sign\,((i)\smallsmile {0}_{i-1}\smallsmile \alpha+(0,1,\ldots,k-1))=(-1)^{i-1} sign\,(\alpha+(0,1,\ldots,k-i-1)).$$
Further, for $\alpha\in G_{k-i},$ 
$${k\choose (i)\smallsmile {0}_{i-1}\smallsmile \alpha}={k\choose (i)\smallsmile \alpha}=
{k\choose i}{k-i\choose \alpha}.$$
Therefore, 
$$\sum_{\alpha\in G_k}sign\,(\alpha+(0,1,\ldots,k-1)){k\choose \alpha}=$$
$$\sum_{i=1}^k 
\sum_{\alpha\in G_{k-i}}(-1)^{i -1} sign\,(\alpha+(0,1,\ldots,k-i-1)){k\choose I}{k-i\choose \alpha}=$$
$$\sum_{i=1}^k (-1)^{i-1} {k\choose i}
\sum_{\alpha\in G_{k-i}}  sign\,(\alpha+(0,1,\ldots,k-i-1)){k-i\choose \alpha}=$$
(by inductive suggestion)
$$\sum_{i=1}^k (-1)^{i-1} {k\choose i}=1.$$

\begin{lm}\label{aaa1}
$$\sum_{i=0}^{l-1}(-1)^i {p\choose i}=(-1)^{l-1}{p-1\choose l-1}.$$
\end{lm}

{\bf Proof.} Induction on $l.$ If $l=1,$ then our statement is evident. Suppose that  it is  true for $l-1\ge 1.$ Then 
$$\sum_{i=0}^{l-1}(-1)^i {p\choose i}=\sum_{i=0}^{l-2}(-1)^i {p\choose i}+(-1)^{l-1} {p\choose l-1}=$$

$$(-1)^{l-1}{p-1\choose l-2}+(-1)^{l-1}{p\choose l-1}=(-1)^{l-1}({p\choose l-1}-{p-1\choose l-2})=$$

$$(-1)^{l-1}{p-1\choose l-1}.$$

\begin{lm}\label{a3} If $k=2l-1,$ then 
$$\sum_{i=1}^l \sum_{\alpha\in G_{l-i}}
sign\,(
(p-l+i)\smallsmile {\bf 0}_{i-1}\smallsmile \alpha\smallsmile {\bf 0}_{l-1}
+(0,1,\ldots,l-1,l,\ldots,2l-2)
) 
{p\choose (p-l+i)\smallsmile \alpha}=$$ $$
{p-1\choose l-1}.$$
If $k=2l,$ then 
$$\sum_{\alpha\in G_{l}}sign\,(
\alpha+(0,1,\ldots,l-1))
{p\choose (p-l)\smallsmile \alpha}=
{p\choose l}.$$
\end{lm}

{\bf Proof.} Let $k=2l-1.$  For  $\alpha\in G_{l-i}$ let 
$\Gamma(\alpha)\in {\bf Z}_0^{2l-1}$ be defined as
$$\Gamma(\alpha)=(p-l+i)\smallsmile {\bf 0}_{i-1}\smallsmile \alpha\smallsmile {\bf 0}_{l-1}
+(0,p-l+1,\ldots,p-1,p+1,\ldots,p+l-1).
$$
Note that 
\begin{equation}\label{ss1}
\Gamma(\alpha)=
(p-l+i,p-l+1,\ldots,p-l+i-1,\alpha_1+p-l+i,\ldots,\alpha_{l-i}+p-1,p+1,\ldots,p+l-1).
\end{equation}

By (\ref{ss1}) 
$$sort(\Gamma(\alpha))=$$
$$
(p-l+1,\ldots,p-l+i-1,p-l+i)\smallsmile sort(\alpha_1+p-l+i,\ldots,\alpha_{l-i}+p-1,p+1,\ldots,p+l-1).
$$
Hence, 
$$sort(\Gamma(\alpha))=\overline{\delta(k)}, \qquad \alpha\in G_{l-i}$$
$$\Updownarrow$$
$$sort(\alpha_1+p-l+i,\ldots,\alpha_{l-i}+p-1,p+1,\ldots,p+l-1)=(p-l+i+1,\ldots,p-1,p,p+1,\ldots,p+l-1).$$
Therefore, the condition $sort(\Gamma(\alpha))=\overline{\delta(k)}$ is equivalent to the condition 
\begin{equation}\label{ss2}
sort(\alpha_1+p-l+i,\ldots,\alpha_{l-i}+p-1)=(p-l+i+1,\ldots,p-1,p).
\end{equation}

By (\ref{ss1})
$$sign\,\Gamma(\alpha)=$$
$$(-1)^{i-1} sign\,(p-l+1,\ldots,p-l+i,
\alpha_1+p-l+i,\ldots,\alpha_{l-i}+p-1,p+1,\ldots,p+l-1).$$
Therefore, by (\ref{ss2}) 
\begin{equation}\label{ss3}
sign\,\Gamma(\alpha)=(-1)^{i-1} sign\,(\alpha_1+p-l+i,\ldots,\alpha_{l-i}+p-1)=
\end{equation}
$$
(-1)^{i-1} sign\,(\alpha_1,\alpha_2+1,\ldots,\alpha_{l-i}+l-i-1).$$
Hence, 
$$\sum_{i=1}^l \sum_{\alpha\in G_{l-i}} sign\, \Gamma(\alpha) 
{p\choose (p-l+i)\smallsmile \alpha}=$$
(by (\ref{ss3}) )
$$\sum_{i=1}^l  \sum_{\alpha\in G_{l-i}}  (-1)^{i-1}sign\, (\alpha+ (0,1,\ldots,l-i-1))
{p\choose l-i}{l-i\choose \alpha}=$$

$$\sum_{i=1}^l  (-1)^{i-1} {p\choose l-i} \sum_{\alpha\in G_{l-i}} sign\, (\alpha+ (0,1,\ldots,l-i-1))
{l-i\choose \alpha}=$$
 (by Lemma \ref{a2} )
$$\sum_{i=1}^l  (-1)^{i-1} {p\choose l-i}=$$
$$\sum_{j=0}^{l-1}(-1)^{l-j-1}{p\choose j}=$$
(by Lemma \ref{aaa1})
$${p-1\choose l-1}.$$
So, our Lemma in case of odd $k$ is proved. 

Let $k=2l.$ Then 
$$\sum_{\alpha\in G_{l}}sign\,(
\alpha+(0,1,\ldots,l-1))
{p\choose (p-l)\smallsmile \alpha}=
$$
$$\sum_{\alpha\in G_{l}}sign\,(
(\alpha_1,\alpha_1+1,\ldots,\alpha_l+l-1)
{p\choose l}{l\choose \alpha}=
$$
$${p\choose l}\sum_{\alpha\in G_{l}}sign\,(
(\alpha_1,\alpha_1+1,\ldots,\alpha_l+l-1){l\choose \alpha}=$$
(by Lemma \ref{a1})
$${p\choose l}.$$
Our Lemma is proved completely. 

\begin{lm} \label{a10} 
Let $\mu_{k}$ be the coefficient at $a^{\delta(k-1)}$ of the element
 $a\der^p (a^{\delta(k-2)}),$ if $k>1,$ and $\mu_1=1.$   If $1\le k\le 2p,$ then   
$$\mu_{k}=
\left\{ \begin{array}{cl}
{p\choose l}, & \mbox{ if $k=2l+1$ is odd,}\\
&\\
{p-1\choose l-1}, & \mbox{ if $k=2l$ is even.}
\end{array}\right. $$
\end{lm}

{\bf Proof.} Follows from Lemmas \ref{a1} and \ref{a3}.

{\bf Example.} If $p=5,$ then 
$$
\begin{array}{|c|l|c|c|}
\hline
k&\delta(k-1)&\mu_{k}\\
\hline
1&(0)&1\\
2&(0,5)&1\\
3&(0,4,6)&5\\
4&(0,4,5,6)&4\\
5&(0,3,4,6,7)&10\\
6&(0,3,4,5,6,7)&6\\
7&(0,2,3,4,6,7,8)&10\\
8&(0,2,3,4,5,6,7,8)&4\\
9&(0,1,2,3,4,6,7,8,9)& 5\\
10&(0,1,2,3,4,5,6,7,8,9)& 1\\
\hline
\end{array}
$$

The following two lemmas can be proved in a similar way as  Lemmas \ref{a1} and \ref{a3}.

\begin{lm} \label{maximal1}
Let $\delta_1(k)$ be maximal element in $E_{k+1,0}(pk-1).$ Then 
$$\delta_1(k)=\left\{\begin{array}{ll}
(0,p-l,\; p-l+2,\ldots,p+l-1), & \mbox{ if $k=2l,$}\\
(0,p-l,p-l+1,\ldots, p, \; p+2,\ldots, p+l), & \mbox{ if $k=2l+1.$}
\end{array}\right.$$
\end{lm}

\begin{lm} \label{delta1}
Let $\gamma_k$ be coefficient at $a^{\delta_1(k-1)}$ of 
$a\der^{p-1}(a^{\delta(k-2)}).$ Then 
$$\gamma_k=p{p-1\choose \lfloor (k-2)/2\rfloor}$$
if $2\le k\le 2p-1.$
\end{lm}

\begin{lm}\label{leader}
Let $\nu_{k}$ be coefficient at $a^{\delta(k-1)}$ of the element $(a\der^p)^{k-1} (a).$  
Then 
 $$leader((a\der^p)^{k})=\nu_k a^{\delta(k-1)}\der^p.$$  
\end{lm}

{\bf Proof.} Follows from Lemma \ref{maximal}.

\begin{lm}  \label{a11} For any $0\le k\le 2p,$ 
$$\nu_k\ge \mu_k\nu_{k-1}.$$   
(Definition of $\mu_k$ see Lemma {\rm \ref{a10}},
and definition of $\nu_k$ see Lemma {\rm \ref{leader}}).
\end{lm}

{\bf Proof.} By  Lemmas \ref{bir1} coefficient at $a^{\delta(k-1)}$ of the element 
$(a\der^p)^{k-1} (a)$ is a non-negative integer that is no less than 
another  non-negative integer 
$(a\der^p)^{k-1} (\nu_{k-1}a^{\delta(k-2)}).$
By Lemma \ref{a10} the last number is equal to  $\nu_{k-1}\mu_k.$

{\bf Example.} Let $p=3.$ Then
$$\mu_1=1,\mu_2=1,\mu_3=3,\mu_4=2,\mu_5=3,\mu_6=1$$
and 
$$(a\der^3)^2=3 a^{(0, 1)}\der^5 + 3 a^{(0, 2)}\der^4 + a^{(0, 3)}\der^3,$$
$$leader((a\der^3)^2)=a^{(0, 3)}\der^3,\quad \nu_2=1,$$

$$(a\der^3)^3=
18 a^{(0, 1, 2)}\der^6 + 27 a^{(0, 1, 3)}\der^5 + 15 a^{(0, 1, 4)}\der^4 + 
 3 a^{(0, 1, 5)}\der^3 + 9 a^{(0, 2, 3)}\der^4 + 3 a^{(0, 2, 4)}\der^3,$$
$$leader((a\der^3)^3)= 3 a^{(0, 2, 4)}\der^3,\quad \nu_3=3,$$

$$(a\der^3)^4=
126 a^{(0, 1, 2, 3)}\der^6 + 189 a^{(0, 1, 2, 4)}\der^5 + 
 99 a^{(0, 1, 2, 5)}\der^4 + $$ $$18 a^{(0, 1, 2, 6)}\der^3 + 
 75 a^{(0, 1, 3, 4)}\der^4 + 24 a^{(0, 1, 3, 5)}\der^3 + 6 a^{(0, 2, 3, 4)}\der^3,$$
$$leader((a\der^3)^4)=6 a^{(0, 2, 3, 4)}\der^3,\quad \nu_4=6,$$

$$(a\der^3)^5=
432 a^{(0, 1, 2, 3, 4)}\der^5 + 432 a^{(0, 1, 2, 3, 5)}\der^4 + 
 108 a^{(0, 1, 2, 3, 6)}\der^3 + 90 a^{(0, 1, 2, 4, 5)}\der^3,$$
$$leader((a\der^3)^5)=90 a^{(0, 1, 2, 4, 5)}\der^3,\quad \nu_5=90,$$

$$(a\der^3)^6= 
90 a^{(0, 1, 2, 3, 4, 5)}\der^3.$$
$$leader((a\der^3)^6)=(a\der^3)^6=a^{(0, 3)}\der^3,\quad \nu_6=90.$$

\begin{lm}\label{s2p} For any $X_1,\ldots,X_N\in A_1^{(p)},$
$$s_N(X_1,\ldots,X_N)=0,$$
if $N>2p$ and 
$$s_{2p}(\der^p,x\der^p,x^2/2\,\der^p,\ldots,x^{2p-1}/(2p-1)!\,\der^p)=\lambda_p\der^p.$$
\end{lm}

{\bf Proof.}  Suppose that 
$X_i=u_i\der^p,$ where $u_i\in K[x].$  Let  us make specialization of $a$ in super-algebra  $\mathcal U.$ Take 
$a=(\sum_{i=1}^N u_i\xi_i)\der^p,$ where $\xi_i$ are odd super-generators. Then 
$$(a\der^p)^N=s_N(u_1\der^p,\ldots,u_N\der^p)\xi_1\cdots \xi_N.$$ 
 By Lemma \ref{N>2p} $(a\der^p)^N=0,$ if $N>2p.$ Therefore, $s_N=0$ is identity if $N>2p.$ 
 
 Now consider the case $N=2p.$ 
 Set $a=\sum_{i=0}^{2p-1} x^{i}/i! \xi_{i+1}$ where $\xi_i$ are odd elements and $\der$ acts on $x^i$ as usual polynomials, $\der(x^i)=ix^{i-1}.$ 
Then 
$$(a\der^p)^{2p}=s_{2p}(\der^p,x\der^p,x^2/2\der^p,\ldots,x^{2p-1}/(2p-1)!\der^p)
\xi_1\xi_2\cdots \xi_{2p}$$
Further,
$$a^{(0,1,2\ldots,2p-1)}=\der^0(a)\der^1(a)\cdots \der^{2p-1}(a)=$$
$$
(\sum_{i=0}^{2p-1} x^{i}/i! \xi_{i+1})
(\sum_{i=0}^{2p-1} x^{i-1}/(i-1)! \xi_{i+1})
\cdots 
(\xi_{2p-1}+x \xi_{2p})
\xi_{2p}=$$
$$
\xi_1\xi_2\cdots\xi_{2p}.$$
Therefore, by Lemma \ref{N>2p} 
$$s_{2p}(\der^p,x\der^p,x^2/2\,\der^p,\ldots,x^{2p-1}/(2p-1)!\der^p)
\xi_1\xi_2\cdots \xi_{2p}=(a\der^p)^{2p}=
\lambda_p \xi_1\xi_2\cdots \xi_{2p}\der^p.$$
Hence
$$s_{2p}(\der^p,x\der^p,x^2/2\,\der^p,\ldots,x^{2p-1}/(2p-1)!\,\der^p)=\lambda_p\der^p.$$

\section{Equivalence of left-commutative and right-commutative identities}

\begin{lm} \label{lcom and rcom}
$(2n-2,1)$-type and $(1,2n-2)$-type identities are equivalent.
\end{lm}

{\bf Proof.} We have to prove that any $n$-algebra $(A,\psi)$ with $(2n-2,1)$-type identity
$$lcom=0,$$
where
$$lcom(t_1,\ldots, t_{2n-1})=\sum_{\sigma\in S^{(2n-2,1)}} sign\,\sigma\, \psi(t_{\sigma(1)},\ldots,
t_{\sigma(n-1)}, \psi(t_{\sigma(n)},\ldots,t_{\sigma(2n-2)},t_{\sigma(2n-1)})),$$
satisfies the identity
$$rcom=0,$$ 
where 
$$rcom(t_1,\ldots,t_{2n-1})=\sum_{\sigma\in S^{(1,2n-2)}} sign\,\sigma\, \psi(t_1,t_{\sigma(2)},\ldots,t_{\sigma(n-1)},\psi(t_{\sigma(n)},\ldots,t_{\sigma(2n-1)})),$$
and vice versa, any  $n$-ary algebra with identity $rcom=0$ satisfies also the identity $lcom=0.$

Let us prove that  
\begin{equation}\label{rtol}
n\, rcom(t_1,\ldots,t_{2n-1})=rcom_1(t_1,\ldots,t_{2n-1}),
\end{equation}

\begin{equation}\label{ltor}
(n-1)\,lcom(t_1,\ldots,t_{2n-1})=lcom_1(t_1,\ldots,t_{2n-1}),
\end{equation}
where
$$rcom_1(t_1,\ldots,t_{2n-1})=$$
$$\sum_{i=2}^{2n-1}(-1)^{i+1}\,  lcom(t_1,\ldots,\hat{t_i},\ldots,t_{2n-1},t_i)- (n-1) \, lcom(t_2,\ldots,t_{2n-1},t_1),$$

$$lcom_1(t_1,\ldots,t_{2n-1})=$$ $$
\sum_{i=1}^{2n-2}(-1)^{i+1}\, rcom(t_i, t_1,\ldots,\hat{t_i},\ldots,t_{2n-1})- (n-2)\, rcom(t_{2n-1},t_1,\ldots,t_{2n-2}).
$$

Note that $rcom(t_1,\ldots,t_{2n-1})$ and $rcom_1(t_1,\ldots,t_{2n-1})$ are skew-symmetric under $2n-2$ variables $t_2,\ldots,t_{2n-1}.$ Therefore, it is enough to prove that coefficients at 
$\psi(t_1,\ldots,t_{n-1},\psi(t_n,\ldots,t_{2n-2},t_{2n-1}))$
and $\psi(t_2,\ldots,t_{n},\psi(t_1,t_{n+1},\ldots,t_{2n-1}))$
  of  $rcom(t_1,\ldots,t_{2n-1})$ and $rcom_1(t_1,\ldots,t_{2n-1})$ are equal. 

It is easy to see that, if $n\le i\le 2n-1,$ then the coefficient at $\psi(t_1,\ldots,t_{n-1},\psi(t_n,\ldots,t_{2n-1})) $ of
 $$(-1)^{i+1}\,  lcom(t_1,\ldots,\hat{t_i},\ldots,t_{2n-1},t_i)$$ is equal to $1.$ If $1\le i<n,$ then  
 this coefficient is $0.$ Therefore, the coefficient at $\psi(t_1,\ldots,t_{n-1},\psi(t_n,\ldots,t_{2n-1}))$ of $rcom_1(t_1,\ldots,t_{2n-1})$ is equal to $n.$ 
 
 Further, if $n\le i\le 2n-1,$ then the coefficient at $\psi(t_2,\ldots,t_{n},\psi(t_1,t_{n+1},\ldots,t_{2n-1})) $ of
 $$(-1)^{i+1}\,  lcom(t_1,\ldots,\hat{t_i},\ldots,t_{2n-1},t_i)$$ is equal to $0.$ If $1\le i<n,$ then  
 this coefficient is $1.$ Therefore, the coefficient at $\psi(t_2,\ldots,t_{n},\psi(t_1,t_{n+1},\ldots,t_{2n-1}))$ of $rcom_1(t_1,\ldots,t_{2n-1})$ is equal to $0.$  
 
 Hence,  relation (\ref{rtol}) is proved completely. 

By similar arguments one establishes (\ref{ltor}).

Relations (\ref{rtol}) and (\ref{ltor}) show that identities $rcom$ and $lcom$ are equivalent.

\section{Proof of Theorem \ref{main}}

By Lemma \ref{s2p} $s_N=0$ is identity on $A_1^{(p)}$ if $N>2p.$
By Lemma \ref{a11}
$$
\lambda_p=\nu_{2p}\ge \mu_{2p}\cdots \mu_{2}\nu_1>0.
$$
Therefore, by Lemma \ref{s2p} $s_{2p}=0$ is not polynomial identity and $s_{2p}$ induces on $A_1^{(p)}$ a non-trivial $2p$-commutator.

By Lemma \ref{a11} for any $1\le k\le 2p-2$
$$\nu_k\ge \mu_{k}\cdots \mu_{2}\nu_1>0.$$
Therefore, by  Lemmas \ref{bir1},  \ref{maximal1} and \ref{delta1} the differential 
$(p+1)$-th order parts of $(a\der^p)^k$  are non-zero for any $2\le k\le 2p-1.$ Therefore, $s_k$ is not well-defined on $A_1^{(p)}.$

Suppose that $A_1^{(p)}$ has identity of degree no more than $2p.$ Then it has skew-symmetric 
multi-linear consequence. In particular, it has a skew-symmetric polynomial identity of degree $2p.$ But $s_{2p}=0,$ as we mentioned above, is not identity. Contradiction. 

Suppose that $I$ is a non-trivial ideal of $A_1^{(p)}$ under $2p$-commutator $s_{2p}.$ Take $0\ne X=u\der^p\in I$ with minimal degree $s=deg\, u.$ Let us prove that $s=0$ and $X=\eta\der^p\in I$ 
for some $0\ne \eta\in K.$ Suppose that it is not true, and $s>0.$ If $s\ge 2p-1,$ then by Lemma \ref{s2p} 
$$s_{2p}(\der^p,x\der^p,\ldots,x^{2p-2}\der^p,X)=\lambda_p{s\choose {2p-1}} \prod_{i=0}^{2p-1}i!\,  x^{s-2p+1}\der^p\in I,$$
or,
$$x^{s-2p+1}\der^p\in I.$$
We obtain contradiction with minimality of $s.$  If $0<s<2p-1,$ then
$$s_{2p}(\der^p,x\der^p,\ldots,x^{s-1}\der^p,X,x^{s+1}\der^p,\ldots,x^{2p-1}\der^p)=\lambda_p 
 \prod_{i=0}^{2p-1}i!\,  \der^p\in I,$$
or,
$$\der^p\in I.$$ Once again we obtain contradiction with minimality of $s.$ 

So, we establish that $X=\eta\der^p\in I,$ for some $0\ne \eta\in K.$ Then for any $l\ge 0,$ 
$$s_{2p}(X,x\der,\ldots,x^{2p-2}\der^p,x^{l+2p-1}\der^p)=\eta \lambda_p{l+2p-1\choose 2p-1} \prod_{i=0}^{2p-1}i!\, 
x^l\der^p\in I.$$
In other words, $x^l\der^p\in I$ for any $l\ge 0.$ This means that $I=A_1^{(p)}.$ 
So, $(A_1^{(p)},s_{2p})$ is simple $2p$-algebra.

By Theorem 1.1 (ii) of \cite{DzhVronskian} the algebra $(A_n(p),s_{2p})$ is left-commutative. 
Presentation of $2p$-commutator  as a Vronskian up to scalar $\lambda_p$ follows from Lemma \ref{s2p}. 


\section{Expressions for  $\lambda_p$}

In this section we give some formulas for  $\lambda_p.$ 
For $s>0$ let us define a polynomial 
$$f_s(x_1,\ldots,x_{2p-1})
=$$
$$\frac{\sum_{\sigma\in Sym_{2p}}sign\,\sigma\,\big(x_{ \sigma(1)}(x_{\sigma(1)}+x_{\sigma(2)})\cdots (x_{\sigma(1)}+x_{\sigma(2)}+\cdots+x_{\sigma(2p-1)})\big)^s}{\prod_{1\le i<j\le 2p}(x_i-x_j).}$$

\bigskip

\noindent {\it  Then $f_s(x_1,\ldots,x_{2p-1})$ is a symmetric polynomial of degree $(2p-1)(s-p).$ In particular, 
$f_p(x_1,\ldots,x_{2p-1})=\lambda_p$ is contsant. The number $\lambda_p$ appears in calculating of $2p$-commutator, }
\bigskip

$$s_{2p}(u_1\der^p,\cdots,u_{2p}\der^p)=\lambda_p
\left|\begin{array}{cccc}u_1&u_2&\cdots&u_{2p}\\
\der(u_1)&\der(u_2)&\cdots&\der(u_{2p})\\
\vdots&\vdots&\cdots&\vdots\\
\der^{2p-1}(u_1)&\der^{2p-1}(u_2)&\cdots&\der^{2p-1}(u_{2p})\\
\end{array}\right|\der^p.$$

\bigskip

\noindent Then 
$$
\lambda_p=\frac{\sum_{\sigma\in Sym_{2p}}sign\,\sigma\,( \sigma(1)(\sigma(1)+\sigma(2))\cdots (\sigma(1)+\sigma(2)+\cdots+\sigma(2p-1)))^p}{\prod_{1\le i<j\le 2p}(i-j).}$$

\noindent For example,
$$\lambda_1=1,\lambda_2=2, \lambda_3=90, \lambda_4=586656, \lambda_5=1915103977500.$$
$$\lambda_6=7886133184567796056800.$$


Another way to calculate $\lambda_p.$ Let ${\mathcal M}_p$ be set of matrices $M=(m_{i,j})$ of order $(2p-1)\times (2p-1)$ such that 
\begin{itemize}
\item $m_{i,j}\in {\bf Z}_0$
\item $m_{i,j}=0 \mbox{ if $i>j$}$
\item sums by rows are constant, $\sum_{j=1}^{2p-1}m_{i,j}=p \mbox{ for any $i$}$
\item sums by columns $r_j=\sum_{i=1}^{2p-1}m_{i,j},$ are positive and different for all $j=1,2,\ldots,2p-1.$ 
\end{itemize}
In particular, 
$$M=(m_{i,j})\in {\mathcal M}_p\Rightarrow m_{1,1}=r_1>0 \mbox{ and }  m_{2p-1,2p-1}=p.$$
For $M\in{\mathcal M}_p$ denote by $r(M)$ the permutation $r_1\ldots r_{2p-1}$ constructed by column sums. 

{\bf Example.} $p=2.$ Then 
$${\mathcal M}_2=\{A=\left(
\begin{array}{ccc}
1&1&0\\0&1&1\\0&0&2\end{array}\right), B=
\left(
\begin{array}{ccc}
1&0&1\\0&2&0\\0&0&2\end{array}\right), C=
\left(
\begin{array}{ccc}
1&1&0\\0&2&0\\0&0&2\end{array}\right), D=
\left(
\begin{array}{ccc}
2&0&0\\0&1&1\\0&0&2\end{array}\right)\}.$$

$$r(A)=123, \;  r(B)=123, \; r(C)=132, \;  r(D)=213. \quad \square$$

\bigskip

{\it  If $M\in {\mathcal M}_p,$ then a sequence $r_1\ldots r_{2p-1}$ induces a permutation, where $r_i=\sum_{j}m_{i,j}$ are sums by columns. 
In particular, $1\le r_i\le 2p-1$ for any $1\le i\le 2p-1.$ Then 
$$\lambda_p=\sum_{M\in {\mathcal M}_p} sign\,r(M) \,\prod_{i=1}^{2p-1} {p\choose m_{i,1} , \ldots, m_{i,2p-1}},$$

$$\lambda_p=
\frac{p!^{2p-1}}{\prod_{j=1}^{2p-1}j!}
\sum_{M \in {\mathcal M}_p} sign\,r(M) \prod_{j} {r_j\choose m_{1,j},\ldots, m_{j,j}}.$$
}

\bigskip

\noindent Here 
$${n\choose n_1,\ldots n_k}=\frac{n!}{n_1!\cdots n_k!}$$
is a multinomial coefficient.


\begin{thebibliography}{10}

\bibitem{Amitzur-Levitzky} A.S. Amitsur,  J. Levitzki, {\em Minimal identities for algebras}, Proc. AMS, {\bf 1}(1050), pp. 449Ð463. 

\bibitem{Dzh} A.S. Dzhumadil'daev,  {\em $N$-commutators,} Comm.Math. Helvetici, {\bf 79}(2004),
No.3, p. 516-553.

\bibitem{DzhVronskian} A.S. Dzhumadil'daev, 
 {\em n-Lie Structures That Are Generated by  Wronskians,}
Sibirskii Matematicheskii Zhurnal, {\bf 46}(2005),
No. 4, pp. 759-773, 2005 =engl. transl.
Siberian Mathematical Journal, {\bf  46}(2005), No.4, pp. 601 - 612

\bibitem{hanlon}  P. Hanlon and M. Wachs, {\em On Lie $k$-algebras,} 
Adv. Math. 113 (1995), 206Ð236.

\bibitem{Martindale} W. S. Martindale,  {\em III Prime rings satisfying a generalized polynomial identity,} J. Algebra  {\bf 12}(1969), 576 - 584.

\end{thebibliography}
\end{document}